\AddToHook{begindocument/before}{} 
\DeclareHookRule{begindocument}{hyperref}{before}{my-mdpi} 

\documentclass[preprints,article,accept,moreauthors,pdftex]{my-mdpi}

\firstpage{1} 
\articlenumber{867}
\doinum{10.3390/axioms12090867} 
\pubvolume{12}
\issuenum{9}
\pubyear{2023}
\copyrightyear{2023}
\externaleditor{Academic Editor: Giuseppe Viglialoro}
\datereceived{10 June 2023} 
\daterevised{4 September 2023} 
\dateaccepted{6 September 2023} 
\datepublished{8 September 2023} 
\hreflink{https://\linebreak doi.org/10.3390/axioms12090867} 

\pdfoutput=1 

\Title{An Analytic Method to Determine the Optimal Time for the Induction Phase of Anesthesia}

\TitleCitation{An Analytic Method to Determine the Optimal Time for the Induction Phase of Anesthesia}

\Author{Mohamed A. Zaitri $^{1,2}$\orcidA{}, 
Cristiana J. Silva $^{1,3,}$*\orcidB{} 
and Delfim F. M. Torres $^{1,4}$\orcidC{}}

\AuthorNames{Mohamed A. Zaitri, Cristiana J. Silva and Delfim F. M. Torres}

\AuthorCitation{Zaitri, M.A.; Silva, C.J.; Torres, D.F.M.}

\address{$^{1}$ \quad Center for Research and Development in Mathematics and Applications (CIDMA),\newline 
Department of Mathematics, University of Aveiro, 3810-193 Aveiro, Portugal;\newline 
zaitri@ua.pt (M.A.Z.); delfim@ua.pt (D.F.M.T.)\\
$^{2}$ \quad Department of Mathematics, University of Djelfa,  Djelfa 17000, Algeria\\
$^{3}$ \quad Department of Mathematics, ISCTE---Instituto Universit\'{a}rio 
de Lisboa, 1649-026 Lisbon, Portugal\\
$^{4}$ \quad Research Center in Exact Sciences (CICE), Faculty of Sciences and Technology (FCT),\newline 
University of Cape Verde (Uni-CV), Praia 7943-010, Cape Verde; delfim@unicv.cv}

\corres{Correspondence: cristiana.joao.silva@iscte-iul.pt}

\abstract{We obtain an analytical solution for the time-optimal 
control problem in the induction phase of anesthesia. Our solution 
is shown to align numerically with the results obtained from 
the conventional shooting method. The induction phase of anesthesia 
relies on a pharmacokinetic/pharmacodynamic (PK/PD) model proposed 
by Bailey and Haddad in 2005 to regulate the infusion of propofol. 
In order to evaluate our approach and compare it with existing 
results in the literature, we examine a minimum-time problem 
for anesthetizing a patient. By applying the Pontryagin minimum 
principle, we introduce the shooting method as a means to solve 
the problem at hand. Additionally, we conducted numerical 
simulations using the MATLAB computing environment.
We solve the time-optimal control problem using our newly 
proposed analytical method and discover that the optimal 
continuous infusion rate of the anesthetic and the minimum 
required time for transition from the awake state to an 
anesthetized state exhibit similarity between the two methods. 
However, the advantage of our new analytic method lies 
in its independence from unknown initial 
conditions for the adjoint variables.}



\begin{document}
	
\section{Introduction} 

Based on Guedel's classification,
the first stage of anesthesia is the induction phase, 
which begins with the initial administration of anesthesia and ends 
with loss of consciousness~\cite{Evers}. Millions of people safely 
receive several types of anesthesia while undergoing medical procedures: 
local anesthesia, regional anesthesia, general anesthesia, and sedation~\cite{Singh}. 
However, there may be some potential complications of anesthesia including anesthetic awareness, 
collapsed lung, malignant hyperthermia, nerve damage, and postoperative delirium. 
Certain factors make it riskier to receive anesthesia, including advanced age, 
diabetes, kidney disease, heart disease, high blood pressure, and smoking~\cite{Merry}.  
To avoid the risk, administering anesthesia should be 
carried out on a scientific basis, based on modern pharmacotherapy, 
which relies on both pharmacokinetic (PK) and pharmacodynamic (PD) 
information~\cite{Beck}. Pharmacokinetics is used to describe the absorption 
and distribution of anesthesia in body fluids, resulting from the administration 
of a certain anesthesia dose. Pharmacodynamics is the study of the effect 
resulting from anesthesia~\cite{Meibohm}. Multiple mathematical models 
were already presented to predict the dynamics of the 
pharmacokinetics/pharmacodynamics (PK/PD) models 
\cite{Absalom,Enlund,MR4502336,MR4482075}. 
Some of these models were implemented following 
different methods~\cite{MR4478523,Singh,MR4539398}.  

The parameters of PK/PD models were fitted by Schnider et al. in~\cite{Schnider}. 
In~\cite{Absalom}, the authors study pharmacokinetic models for propofol,  
comparing Schnider et al. and Marsh et al. models~\cite{Marsh}. 
The authors of~\cite{Absalom} conclude that Schnider's model should always be used 
in effect-site targeting mode, in which larger initial doses are administered 
but smaller than those obtained from Marsh's model. However, users of the Schnider model 
should be aware that in the morbidly obese, the lean body mass (LBM) equation can 
generate paradoxical values, resulting in excessive increases in maintenance infusion rates 
\cite{Schnider}. In~\cite{Zabi}, a new strategy is presented to develop a robust control 
of anesthesia for the maintenance phase, \textls[-15]{taking into account the saturation of the actuator. 
The authors of~\cite{Said} address the problem of}\linebreak  \textls[-35]{optimal control of the induction phase. 
For other related works, see~\cite{MR4410555,MR4502336} and references therein.}

Here, we consider the problem proposed in~\cite{Said}, 
to transfer a patient from a state consciousness to unconsciousness. 
We apply the shooting method \cite{Bock}
using the Pontryagin minimum principle~\cite{Pontryagin},
correcting some inconsistencies found in~\cite{Said} related
with the stop criteria of the algorithm and the numerical computation
of the equilibrium point. Secondly, we provide a new different analytical method
to the time-optimal control problem for the induction phase of anesthesia. 
While the shooting method, popularized by Zabi et al.~\cite{Said}, 
is widely employed for solving such control problems and determining the minimum time, 
its reliance on Newton's method makes it sensitive to initial conditions. The shooting method's 
convergence is heavily dependent on the careful selection of initial values, particularly 
for the adjoint vectors. To overcome this limitation, we propose an alternative approach, 
which eliminates the need for initial value selection and convergence analysis. Our method offers 
a solution to the time-optimal control problem for the induction phase of anesthesia, 
free from the drawbacks associated with the shooting method. Furthermore, we propose that 
our method can be extended to other PK/PD models to determine optimal timings for drug administration. 
To compare the methods, we perform numerical simulations to compute the minimum time 
to anesthetize a man of 53 years, 77~kg, 
and 177~cm, as considered in~\cite{Said}.
We find the optimal continuous infusion rate of the anesthetic and the minimum 
time that needs to be chosen for treatment, showing that both 
the shooting method of~\cite{Said} and the one proposed here coincide. 

This paper is organized as follows. In Section~\ref{Section:PK:PD:model},
we recall the pharmacokinetic and pharmacodynamic model of Bailey and Haddad 
\cite{Bailey}, the Schnider model~\cite{Schnider}, the bispectral index (BIS), 
and the equilibrium point~\cite{Zabi}. Then, in Section~\ref{Minimum:time:problem}, 
a time-optimal control problem for the induction phase of anesthesia is posed
and solved both by the shooting and analytical methods. 
Finally, in Section~\ref{Numerical:example}, we compute the parameters of the model 
using the Schnider model~\cite{Schnider}, and we illustrate the results of the time-optimal 
control problem through numerical simulations. We conclude that the optimal continuous infusion 
rate for anesthesia and the minimum time that should be chosen for this treatment
can be found by both shooting and analytical methods. The advantage of the new method
proposed here is that it does not depend on the concrete initial conditions, while
the shooting method is very sensitive to the choice of the initial conditions
of the state and adjoint variables. We end with Section~\ref{sec:conc} of conclusions,
pointing also some directions for future research.


\section{The PK/PD Model}
\label{Section:PK:PD:model}

\begin{figure}[H]
\vspace{-4pt}

\includegraphics[scale=0.75]{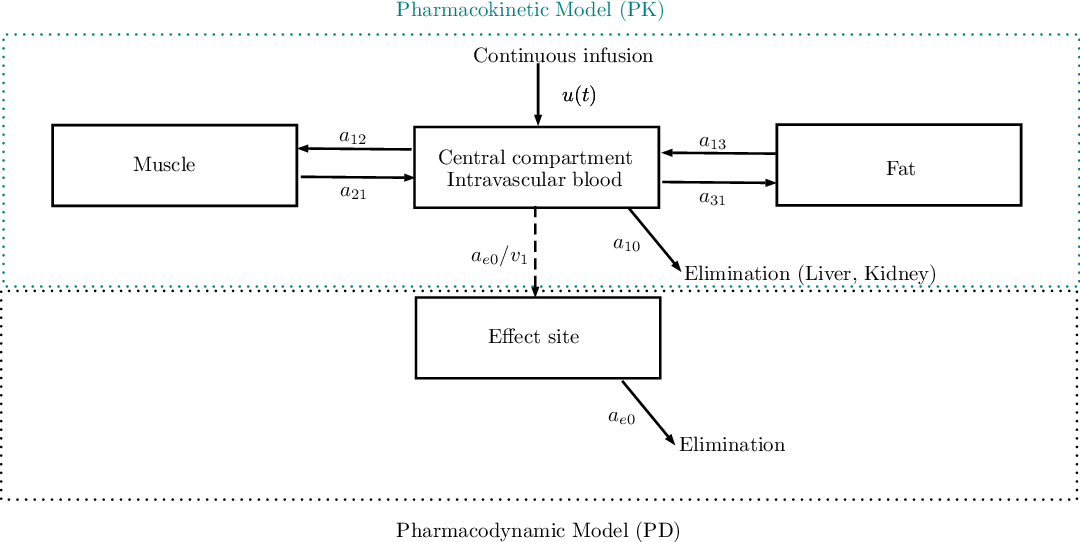}
\caption{Schematic diagram of the PK/PD model with the effect site compartment 
of Bailey and Haddad~\cite{Bailey}.}
\label{schema01}
 
\end{figure}

The pharmacokinetic/pharmacodynamic (PK/PD) model consists of four compartments: 
intravascular blood $(x_1(t))$, muscle $(x_2(t))$, fat $(x_3(t))$, and effect site $(x_4(t))$.  
The effect site compartment (brain) is introduced to account for the finite equilibration time 
between central compartment  and central nervous system concentrations~\cite{Bailey}. 
This  model is used to describe  the circulation of drugs in a patient's body, being expressed
by a four-dimensional dynamical system as follows:
\begin{equation} 
\label{model:PK/PD} 
\left\{
\begin{array}{l l}
\dot{x}_1(t)= -(a_{1\,0}+a_{1\,2}+a_{1\,3})\,x_1(t)+a_{2\,1}\,x_2(t)+ a_{3\,1}\,x_3(t)+u(t),\\		
\dot{x}_2(t)= a_{1\,2}\,x_1(t) -a_{2\,1} \, x_2(t),\\	
\dot{x}_3(t)= a_{1\,3}\,	x_1(t)-a_{3\,1}\,x_3(t),\\		
\dot{x}_4(t)= \frac{a_{e\,0}}{v_1}\,x_1(t) -a_{e\,0}\, x_4(t).
\end{array}\right.
\end{equation}
The state variables for system \eqref{model:PK/PD} 
are subject to the following initial conditions:
\begin{equation}
\label{initial:state}
x(0)=\left(x_1(0), x_2(0),x_3(0),x_4(0)\right)= \left(0,0,0,0 \right),
\end{equation}
where $x_1(t),x_2(t),x_3(t)$, and $x_4(t)$ represent, respectively, 
the masses of the propofol in the compartments of blood, muscle, fat, 
and effect site at time $t$. 
The control $u(t)$ is the continuous infusion rate of the anesthetic. The parameters 
$a_{1\,0}$ and $a_{e\,0}$ represent, respectively, the rate of clearance from the central 
compartment and the effect site. The parameters $a_{1\,2}$, $a_{1\,3}$, $a_{2\,1}$, 
$a_{3\,1}$, and $a_{e\,0}/v_1$ are the transfer rates of the drug between compartments. 
A schematic diagram of the dynamical control system \eqref{model:PK/PD} 
is given in Figure~\ref{schema01}. 


\subsection{Schnider's Model}

Following Schnider et al.~\cite{Schnider}, 
the lean body mass (LBM) is calculated using the James
formula, which performs satisfactorily in normal and
moderately obese patients, but not so well for
severely obese cases~\cite{James}. 
The James formula calculates LBM as follows:
\begin{eqnarray}
\text{for Male},\, 
\text{LBM}&=&1.1\times \text{weight}
-128\times\left( \dfrac{\text{weight}}{\text{height}}\right)^2,\\	
\text{for Female},\, \text{LBM}&=&1.07\times \text{weight}
-148\times\left( \dfrac{\text{weight}}{\text{height}}\right)^2.
\end{eqnarray}
The parameters of the PK/PD model \eqref{model:PK/PD} 
are then estimated according to Table~\ref{Table:Schnider}. 
\begin{table}[H]\footnotesize 
\caption{Parameter values for model \eqref{model:PK/PD}
according to Schnider's model~\cite{Schnider}.}
\label{Table:Schnider}
 \tabcolsep=0.823cm
\begin{tabular}{cc} \toprule
\textbf{Parameter} 
& \textbf{Estimation} \\ \midrule
$a_{10}\,(\text{min}^{-1})$ 
& $0.443+0.0107\, (\text{weight}-77)-0.0159\, 
(\text{LBM}-59)+0.0062\, (\text{height}-177)$\\ \midrule
$a_{12}\,(\text{min}^{-1})$ 
& $0.302-0.0056 \, (\text{age}-53)$\\ \midrule
$a_{13}\,(\text{min}^{-1})$ 
& 0.196 \\ \midrule
$a_{21}\,(\text{min}^{-1})$ 
& $\left( 1.29-0.024\, (\text{age}-53)\right) 
/\left( 18.9-0.391\, (\text{age}-53)\right) $ \\ \midrule
$a_{31}\,(\text{min}^{-1})$ & 0.0035 \\ \midrule
$a_{e0}\,(\text{min}^{-1})$ & 0.456 \\ \midrule
$v_1\,(\text{L})$ & 4.27\\ \bottomrule
\end{tabular}
\end{table}


\subsection{The Bispectral Index (BIS)}

The BIS is the depth of anesthesia indicator, which is a signal derived from the EEG analysis 
and directly related to the effect site concentration of $x_4(t)$. It quantifies the level 
of consciousness of a patient from 0 (no cerebral activity) to 100 (fully awake patient), 
and can be described empirically by a decreasing sigmoid function~\cite{Bailey}:
\begin{equation}
BIS(x_4(t))
=BIS_0\left(1-\dfrac{x_4(t)^\gamma}{x_4(t)^\gamma+EC_{50}^\gamma} \right), 
\end{equation}
where $BIS_0$ is the $BIS$ value of an awake patient typically set to 100, 
$EC_{50}$ corresponds to the drug concentration associated
with $50\%$ of the maximum effect, and $\gamma$ is a parameter
modeling the degree of nonlinearity. According to~\cite{Haddad},
typical values for these parameters are $EC_{50}$ = 3.4\,mg/L and $\gamma=3$.


\subsection{The Equilibrium Point}

Following~\cite{Zabi}, the equilibrium point is obtained by
equating the right-hand side of \eqref{model:PK/PD} to zero,
\begin{equation} 
\left\{
\begin{array}{l l}
0 = -(a_{1\,0}+a_{1\,2}+a_{1\,3})\,x_1+a_{2\,1}\,x_2+ a_{3\,1}\,x_3+u, \\		
0 =	a_{1\,2}\,x_1 -a_{2\,1} \, x_2,\\		
0 =	a_{1\,3}\,	x_1-a_{3\,1}\,x_3,\\		
0 =\frac{a_{e\,0}}{v_1}\,x_1 -a_{e\,0}\, x_4,\\
\end{array}\right.
\end{equation}
with the condition
\begin{equation}
x_4=EC_{50}.
\end{equation} 
It results that the equilibrium point 
$x_e=\left(x_{e\,1},x_{e \, 2},x_{e\,  3},x_{e\,  4} \right)$ is given by
\begin{equation}
x_{e\,  1  }=v_1\,EC_{50}, \quad
x_{e\,  2}= \frac{a_{1\,2}\,v_1\,EC_{50}}{a_{2\,1}}, \quad
x_{e\,  3}=\frac{a_{1\,3}\,v_1\,EC_{50}}{a_{3\,1}}, \quad
x_{e\,  4}=EC_{50},	
\end{equation}
and the value of the continuous infusion rate for this equilibrium is
\begin{equation}
u_{e}=a_{1\,0}\,v_1 \,EC_{50}. 
\end{equation}

The fast state is defined by
\begin{equation}
\label{fast:final:state}
x_{eF}(t)=(x_1(t),x_4(t)).
\end{equation}
The control of the fast dynamics is crucial because the BIS is 
a direct function of the concentration at the effect site.  


\section{Time-Optimal Control Problem}
\label{Minimum:time:problem}

Let $x(t)=(x_1(t),x_2(t),x_3(t),x_4(t))\in \mathbb{R}^4$.
We can write the dynamical system \eqref{model:PK/PD}  
in a matrix form as follows:
\begin{equation}
\label{dot:x}
\dot{x}(t)=A\,x(t)+B\,u(t),
\end{equation}
where
\begin{equation}
\label{model:1:matrix}
A=\left(
\begin{array}{cccc}
-(a_{1\,0}+a_{1\,2}+a_{1\,3}) & a_{2\,1} & a_{3\,1} & 0 \\
a_{1\,2} & -a_{2\,1} & 0 & 0 \\
a_{1\,3} & 0 & -a_{3\,1} & 0 \\
\frac{a_{e\,0}}{v_1} & 0 & 0 & -a_{e\,0} 
\end{array}
\right)
\quad \text{ and } \quad
B=\left(
\begin{array}{c}
1 \\
0  \\
0  \\
0 
\end{array}
\right). 
\end{equation}
Here, the continuous infusion rate $u(t)$ is to be chosen so as to transfer the system 
\eqref{model:PK/PD} from the initial state (wake state) to the fast final state 
(anesthetized state) in the shortest possible time. Mathematically, 
we have the following time-optimal control problem~\cite{Said}:
\begin{equation}
\label{Time:Optimal}
\begin{cases}
\min\limits_{t_f} J=\int\limits_0^{t_f}dt,\\
\dot{x}(t)=A\,x(t)+B\,u(t),	
\quad x(0)=(0,0,0,0),\\
C\,x_{eF}(t_f)=x_{eF},\\
0\leq u(t)\leq U_{max},\, 
\quad t\in [0,t_f],\, 
\quad t_f \, \text{is free,}
\end{cases}
\end{equation}
where $t_f$ is the first instant of time that the desired state is reached, 
and $C$ and $x_{eF}$ are given by
\begin{equation}
C=\left(
\begin{array}{cc}
1  & 0 \\
0  & 1
\end{array}
\right), \,\,\,\, 
x_{eF}=(x_{e1},\,x_{e4}),
\end{equation} 
with
\begin{equation}
\label{x:e:f}		
x_{eF}(t_f)=(x_1(t_f),x_2(t_f)).
\end{equation}


\subsection{Pontryagin Minimum Principle} 
\label{subsec:PMP} 

According to the Pontryagin minimum principle (PMP)~\cite{Pontryagin}, 
if $\tilde{u}\in L^1$ is optimal for Problem \eqref{Time:Optimal} 
and the final time $t_f$ is free, then there exists 
$\psi(t) = (\psi_1(t), \ldots , \psi_4(t))$, $t\in [0,t_f]$,
$\psi\in AC([0,t_f];\mathbb{R}^4)$,
called the adjoint vector, such that  
\begin{equation}
\label{adj:system}
\begin{cases}
\dot{x}=\displaystyle \frac{\partial H}{\partial \psi},\\[0.3cm]
\dot{\psi}=-\displaystyle \frac{\partial H}{\partial x},
\end{cases}
\end{equation}
where the Hamiltonian $H$ is defined by
\begin{equation}
\label{def:H}
H(t,x,u,\psi)=1+\psi^T\,(A\,x+B\, u).
\end{equation}
Moreover, the minimality condition
\begin{equation}
\label{Principe:Max:Hamilton}
H(t,\tilde{x}(t),\tilde{u}(t),\tilde{\psi}(t))
=\min\limits_{0\leq u\leq U_{max}}H(t,\tilde{x}(t),u,\tilde{\psi}(t))
\end{equation}
holds almost everywhere on $t\in [0,t_f]$.

Since the final time $t_f$ is free, 
according to the transversality condition of PMP, we obtain:
\begin{equation}
\label{condition:transversality}
H(t_f,x(t_f),u(t_f),\psi(t_f) )=0.
\end{equation}

Solving the minimality condition \eqref{Principe:Max:Hamilton} 
on the interior of the set of admissible controls gives
the necessary condition
\begin{equation}
\label{Optimal:control:u}
\tilde{u}(t)=\begin{cases}
0 & \hbox{if}\,\, \tilde{\psi}_1(t)>0, \\
U_{max}& \hbox{if}\,\, \tilde{\psi}_1(t)<0,
\end{cases}	
\end{equation}
where $\tilde{\psi}_1(t)$ is obtained from the adjoint system
\eqref{adj:system}, that is, $\tilde{\psi}'(t) = -A^T \tilde{\psi}(t)$,
and the transversality condition \eqref{condition:transversality}.
This is discussed in Sections~\ref{Shooting:method} and \ref{Analytical:solution}.


\subsection{Shooting Method}
\label{Shooting:method}	

The shooting method is a numerical technique used to solve boundary 
value problems, specifically in the realm of differential equations and optimal control. 
It transforms the problem into an initial value problem by estimating the unknown 
boundary conditions. Through iterative adjustments to these estimates, the boundary 
conditions are gradually satisfied. In~\cite{Bock}, the authors propose an algorithm 
that addresses numerical solutions for parameterized optimal control problems. This 
algorithm incorporates multiple shooting and recursive quadratic programming, introducing 
a condensing algorithm for linearly constrained quadratic subproblems and high-rank 
update procedures. The algorithm's implementation leads to significant improvements 
in convergence behavior, computing time, and storage requirements. For more on numerical 
approaches to solve optimal control problems, we refer the reader 
to~\cite{Zaitri2019} and references therein.

Using \eqref{adj:system}, \eqref{def:H}, \eqref{condition:transversality},
and \eqref{Optimal:control:u}, we consider the following problem:
\begin{equation}
\label{two:points:boundary}
\begin{cases}
\dot{x}(t)=A\,x(t)+B\,\times\max\,(0,-U_{max}\,sign (\psi_1(t))),	\\
\dot{\psi}(t)=-A^T\,\psi(t),	\\
x(0)=(0,0,0,0),\, x_{1}(t_f)=x_{e1},\, x_{4}(t_f)=x_{e4},\\
\psi(0)\, \text{is free},\, 
H(t_f,x(t_f),\max\,(0,-U_{max}\,sign (\psi_1(t_f))),\psi(t_f) )=0.
\end{cases}
\end{equation}

Let $z(t)=(x(t),\psi(t))$. Then, we obtain the following 
two points' boundary value problem:
\begin{equation}
\label{Z}
\begin{cases}
\dot{z}(t)=A^*z(t)+B^*,\\
\textit{R}(z(0),z(t_f))=0,
\end{cases}
\end{equation}	
where $A^*\in \textit{M}_{8\times 8}(\mathbb{R})$ 
is the matrix given by
\begin{equation}
A^*=\left(
\begin{array}{cc}
A & 0_{4\times 4} \\
0_{4\times 4} & -A^T
\end{array}
\right),
\end{equation} 
$B^*\in \mathbb{R}^8$ is the vector given by
\begin{equation}
B^*=\begin{cases}
\left( 0,\,0,\,0,\,0,\,0,\,0,\,0,\,0\right) & \hbox{if}\,\, \psi_1(t)>0, \\
\left( U_{max},\,0,\,0,\,0,\,0,\,0,\,0,\,0\right) & \hbox{if}\,\, \psi_1(t)<0,
\end{cases}	
\end{equation}
and $R(z(0),z(t_f))$ is given by 
\eqref{initial:state}, \eqref{x:e:f}, 
and \eqref{condition:transversality}.
We consider the following Cauchy problem:
\begin{equation}
\label{Cauchy:problem}
\begin{cases}
\dot{z}(t)=A^*z(t)+B^*,\\
z(0)=z_0.
\end{cases}
\end{equation}	
If we define the shooting function 
$\textit{S}:\, \mathbb{R}^4\longrightarrow\mathbb{R}^3$ by
\begin{equation}
\textit{S}(z_0)=\textit{R}(t_f,z(t_f,z_0)),
\end{equation}
where $z(t,z_0)$ is the solution of the Cauchy problem \eqref{Cauchy:problem},
then the two points' boundary value problem \eqref{two:points:boundary} 
is equivalent to 
\begin{equation}
\label{eq:NM}
\textit{S}(z_0)=0.
\end{equation}
To solve \eqref{eq:NM}, we use Newton's method~\cite{Newton}.


\subsection{Analytical Method}
\label{Analytical:solution}

We now propose a different method to choose the optimal control. 
If the pair $(A, B)$ satisfies the Kalman condition and all 
eigenvalues of matrix $A\in n\times n$ are real, then any extremal control 
has at most $n-1$ commutations on $\mathbb{R}^+$ (at most $n-1$ switching times). 
We consider the following eight possible strategies:

\medskip

\noindent Strategy~1 (zero switching times): 
\begin{equation}
u(t)=U_{max},\, \forall t\in[0,t_f].	
\end{equation}
Strategy 2 (zero switching times):
\begin{equation}
u(t)=0,\, \forall t\in[0,t_f].	
\end{equation}
Strategy 3 (one switching time):
\begin{equation}
\label{Optimal:control:u:2}
u(t)=
\begin{cases}
U_{max}& \hbox{if}\,\, 0\leq t<t_c, \\
0 & \hbox{if}\,\, t_c<t\leq t_f,
\end{cases}	
\end{equation}
where $t_c$ is a switching time.\\[0.3cm]
Strategy 4 (one switching time):
\begin{equation}
u(t)=\begin{cases}
0& \hbox{if}\,\, 0\leq t<t_c, \\
U_{max} & \hbox{if}\,\, t_c<t\leq t_f.
\end{cases}	
\end{equation} 
Strategy 5 (two switching times):
\begin{equation}
u(t)=
\begin{cases}
U_{max} & \hbox{if}\,\, 0<t<t_{c1}, \\
0& \hbox{if}\,\, t_{c1}<t<t_{c2}.\\
U_{max} & \hbox{if}\,\, t_{c2}< t\leq t_f, \\
\end{cases}	
\end{equation} 
where $t_{c1}$ and $t_{c2}$ represent two switching times.\\[0.3cm]
Strategy 6 (two switching times):
\begin{equation}
u(t)
=\begin{cases}
0 & \hbox{if}\,\, 0<t<t_{c1}, \\
U_{max}& \hbox{if}\,\, t_{c1}<t<t_{c2}.\\
0 & \hbox{if}\,\, t_{c2}< t\leq t_f. \\
\end{cases}	
\end{equation} 
Strategy 7 (three switching times):
\begin{equation}
u(t)=
\begin{cases}
U_{max} & \hbox{if}\,\, 0<t<t_{c1}, \\
0& \hbox{if}\,\, t_{c1}<t<t_{c2}.\\
U_{max} & \hbox{if}\,\, t_{c2}< t\leq t_{c3}. \\
0& \hbox{if}\,\, t_{c3}<t<t_{f},\\
\end{cases}	
\end{equation} 
where $t_{c1}$, $t_{c2}$, and $t_{c3}$ represent three switching times.\\[0.3cm]
Strategy 8 (three switching times):
\begin{equation}
u(t)
=\begin{cases}
0 & \hbox{if}\,\, 0<t<t_{c1}, \\
U_{max}& \hbox{if}\,\, t_{c1}<t<t_{c2}.\\
0 & \hbox{if}\,\, t_{c2}< t\leq t_{c3}. \\
U_{max}& \hbox{if}\,\, t_{c3}<t<t_{f}.\\
\end{cases}	
\end{equation} 
Let $x(t)$ be the trajectory associated with the 
control $u(t)$, given by the relation
\begin{equation}
x(t)=\exp(A\,t)\,x(0)+\int\limits_0^{t}\exp(A(t-s))Bu(t)ds,
\end{equation}
where $\exp(A)$ is the exponential matrix of $A$. 

To calculate the switching times $t_c$, $t_{c1}$, $t_{c2}$  
and the final time $t_f$, we have to solve the 
following nonlinear equation:
\begin{equation}
\label{eq:xEF}
\tilde{x}_{eF}(t_f)=(x_{e1},\,x_{e4}).
\end{equation}
We also solve \eqref{eq:xEF} using the Newton method~\cite{Newton}.


\section{Numerical Example}
\label{Numerical:example}

In this section, we use the shooting and analytical methods 
to calculate the minimum time $t_f$ to anesthetize 
a man of 53 years, $77$~kg, and $177$~cm.

The equilibrium point and the flow rate
corresponding to a BIS of 50 are:
\begin{equation}
x_e=(14.518\, \text{mg}  ,\, 64.2371 \, \text{mg}  ,\,  
813.008\, \text{mg}  ,\, 3.4\, \text{mg}),\,\,\, u_e=6.0907\, \text{mg/min}.
\end{equation}

Following the Schnider model, the matrix $A$ 
of the dynamic system \eqref{dot:x} is given by: 
\begin{equation}
A=\left(
\begin{array}{cccc}
 -0.9175  &  0.0683 &   0.0035&         0\\
0.3020 &  -0.0683  &       0    &     0\\
0.1960  &       0 &  -0.0035     &    0\\
0.1068  &       0 &         0  & -0.4560
\end{array}
\right)
\quad \text{ and }\quad 
B=\left(
\begin{array}{c}
1 \\
0  \\
0  \\
0 
\end{array}
\right). 
\end{equation}	

We are interested to solve the following minimum-time control problem:
\begin{equation}
\label{Example1:Time:Optimal}
\begin{cases}
\min\limits_{t_f} J=\int\limits_0^{t_f}dt,\\
\dot{x}(t)=A\,x(t)+B\,u(t),	\quad x(0)=(0,\,0,\,0,\,0),\\
x_{e1}(t_f)= 14.518 \, \text{mg},\quad x_{e4}(t_f)=3.4\, \text{mg},\\
0\leq u(t)\leq 106.0907 ,\quad t\in [0,t_f],
\quad t_f \, \text{is free.}
\end{cases}
\end{equation}


\subsection{Numerical Resolution by the Shooting Method}
\label{sec:SM}

Let $z(t)=(x(t),\psi(t))$. We consider the following Cauchy problem:
\begin{equation}
\begin{cases}
\dot{z}(t)=A^*z(t)+B^*,\\
z(0)=z_0=(0,\,0,\,0,\,0,\,\psi_{01},\,
\psi_{02},\,\psi_{03},\,\psi_{04}),
\end{cases}
\end{equation}	
where
\begin{equation}
A^*=10^{-4}\left(
\begin{array}{cccccccc}
-9175  &  683  & 35 & 0  & 0& 0& 0& 0 \\
3020 &  -683 & 0 & 0   & 0& 0& 0& 0      \\
196  &       0 & -35   & 0& 0& 0& 0  & 0      \\
1068 & 0 & 0 &-456     & 0& 0& 0& 0 \\
0& 0& 0& 0  & 9175  & -3020 & -196  &-1068\\
0& 0& 0& 0  & -683 & 683  &  0       &0\\
0& 0& 0& 0  &-35  & 0       & 35   &0\\
0& 0& 0& 0  &0        & 0       &  0       &456
\end{array}
\right), \,   
\end{equation}
\begin{equation}	
B^*=\left(
\begin{array}{c}
\max\,(0, -106.0907 \,sign(\psi_1(t)))\\
0\\
0\\
0\\
0 \\
0 \\
0 \\
0
\end{array}
\right).
\end{equation}	

The shooting function \textit{S} is given by 
\begin{equation}
\textit{S}(z_0)=(\textit{S}_1(z_0),\,\textit{S}_2(z_0),\,\textit{S}_3(z_0)),
\end{equation}
where
\begin{eqnarray*}
\textit{S}_1(z_0)&=&x_{e1}(t_f)-14.518,\\
\textit{S}_2(z_0)&=&x_{e4}(t_f)-3.4,\\
\textit{S}_3(z_0)&=& 1+\psi^T(t_f)\left(Ax(t_f)
+B\max\,(0,-106.0907 \, sing\, \psi_1(t_f))\right).
\end{eqnarray*}

All computations were performed with the MATLAB numeric
computing environment, version R2020b, using the medium-order method 
and the function \texttt{ode45} (Runge--Kutta method) in order to solve 
the nonstiff differential system \eqref{Z}. 
We have used the variable order method and the function \texttt{ode113} 
(Adams--Bashforth--Moulton method) in order to solve the nonstiff differential system
\eqref{Cauchy:problem}, and the function \texttt{fsolve} in order to solve 
equation $\textit{S}(z_0)=0$. Thus, we obtain that the minimum time is equal to 
\begin{equation}
t_f=1.8397\, \text{min}, 
\end{equation}
with 
\begin{equation}
\label{eq:ic:psi}
\psi^T(0)=( -0.0076,\, 0.0031,\, -0.0393,\,  -0.0374).
\end{equation}


\subsection{Numerical Resolution by the Analytical Method}

The pair $(A,B)$ satisfies the Kalman condition, and the matrix $A$ 
has four real eigenvalues. Then, the extremal control $u(t)$ has at most three 
commutations on $\mathbb{R}^+$. Therefore, let us test the eight strategies 
provided in Section~\ref{Analytical:solution}.
	
Note that the anesthesiologist begins with a bolus injection 
to transfer the patient state from the consciousness state $x(0)$ 
to the unconsciousness state 
$$
x_{eF}=(14.518 ,\,3.4),
$$ 
that is,
\begin{equation}
u(0)=U_{max}= 106.0907\, \text{mg/min}.
\end{equation}
Thus, Strategies 2, 4, 6, and 8 are not feasible here. 
Therefore, in the sequel, we investigate Strategies 1, 3, 5, and 7 only.\\[0.3cm]
Strategy~1: Let $u(t)= 106.0907\, \text{mg/min}$ for all $t\in[0,t_f]$. 
The trajectory $x(t)$, associated with this
control $u(t)$, is given by the following relation:
\begin{equation}
x(t)=\int\limits_0^{t}\exp(A(t-s))BU_{max}ds, \,\, \forall t\in[0,t_f],
\end{equation}
where 
\begin{eqnarray}
\exp(A\,(t-s))	&=&V\, D(t-s)\, V^{-1}
\end{eqnarray}
with
\begin{equation} 
 V=\left(
\begin{array}{cccc}
0&    0.9085&    0.0720&   -0.0058\\
0&   -0.3141&    0.9377&   -0.0266\\
0&   -0.1898&   -0.3395&   -0.9996\\
1&   -0.1997&    0.0187&   -0.0014
\end{array}
\right)
\end{equation}
and
\begin{equation}
D(\tau)=\left(
\begin{array}{cccc}
\exp^{-0.4560\,\tau}   &      0   &      0   &      0\\
0                       &  \exp^{-0.9419\,\tau}  &       0    &     0\\
0   &      0 &  \exp^{-0.0451\,\tau}    &     0\\
0   &      0  &       0  & \exp^{-0.0024 \,\tau}
\end{array}
\right).
\end{equation}  
System \eqref{eq:xEF} takes the form
\begin{equation}
\label{non:linear:equation}
\begin{cases}
x_1(t_f)= 14.518,\\
x_4(t_f)= 3.4,\\
\end{cases}
\end{equation}
and has no solutions. Thus, Strategy 1 is not feasible.\\[0.3cm]
Strategy 3:  Let
$u(t),\,t\in [0,t_f],$ be the control defined by 
\begin{equation}
\label{Optimal:control:u:3}
u(t)=\begin{cases}
106.0907\, \text{mg/min}& \hbox{if}\,\, 0\leq t<t_c, \\
0 & \hbox{if}\,\, t_c<t\leq t_f.
\end{cases}	
\end{equation}
The trajectory $x(t)$ associated with this
control $u(t)$ is given by 
\begin{equation}
x(t)=\begin{cases}
\int\limits_0^{t}\exp(A(t-s))BU_{max}ds 
& \hbox{if}\,\, 0\leq t\leq t_c, \\
\exp(A\,(t-t_c))\,x(t_c)
& \hbox{if}\,\, t_c<t\leq t_f,	
\end{cases}
\end{equation}
where 
\begin{eqnarray}
\exp(A\,(t-t_c))	&=&V\, D(t-t_c)\, V^{-1}.
\end{eqnarray}
To calculate the switching time $t_c$ and the final time $t_f$, 
we have to solve the nonlinear system \eqref{non:linear:equation} 
with the new condition
\begin{equation}
t_c< t_f.
\end{equation}
Similarly to Section~\ref{sec:SM}, all numerical computations 
were performed with MATLAB R2020b using the command \texttt{solve} 
to solve Equation \eqref{non:linear:equation}. 
The obtained minimum time is equal to
\begin{equation}
t_f=1.8397\, \text{min}, 
\end{equation}
with the switching time
\begin{equation}
t_c=0.5467\, \text{min}.
\end{equation}
Strategy 5: Let $u(t)$, $t\in [0,t_f]$, 
be the control defined by the relation
\begin{equation}
\label{control:57}
u(t)=\begin{cases}
106.0907\,\text{mg/min} & \hbox{if}\,\, 0\leq t<t_{c1}, \\
0& \hbox{if}\,\, t_{c1}<t<t_{c2}.\\
106.0907\,\text{mg/min} & \hbox{if}\,\, t_{c2}< t\leq t_f, \\
\end{cases}	
\end{equation} 
where $t_{c1}$ and $t_{c2}$ are the two switching times. 
The trajectory $x(t)$ associated with control 
\eqref{control:57} is given by 
\begin{equation}
x(t)=\begin{cases}
\int\limits_0^{t}\exp(A(t-s))BU_{max}ds & \hbox{if}\,\, 0\leq t\leq t_{c1}, \\
\exp(A\,(t-t_{c1}))\,x(t_{c1})& \hbox{if}\,\, t_{c1}<t\leq t_{c2},	\\
\exp(A\,(t-t_{c2}))\,x(t_{c2})+\int\limits_{t_{c2}}^{t}\exp(A(t-s))BU_{max}ds
& \hbox{if}\,\, t_{c2}<t\leq t_{f}.
\end{cases}
\end{equation}
To compute the two switching times $t_{c1}$ and $t_{c2}$ 
and the final time $t_f$, we have to solve the 
nonlinear system \eqref{non:linear:equation} with
\begin{equation}
\label{switching:time:tc1:tc2}
0\leq t_{c1}\leq t_{c2}\leq t_f.
\end{equation}
It turns out that System \eqref{non:linear:equation} subject to 
Condition \eqref{switching:time:tc1:tc2} has no solution. 
Thus, Strategy~5 is also not feasible.\\[0.3cm]
Strategy 7: Let $u(t)$, $t\in [0,t_f]$, 
be the control defined by the relation
\begin{equation}
\label{control:S7}
u(t)=\begin{cases}
106.0907\,\text{mg/min} & \hbox{if}\,\, 0\leq t<t_{c1}, \\
0& \hbox{if}\,\, t_{c1}<t<t_{c2}.\\
106.0907\,\text{mg/min} & \hbox{if}\,\, t_{c2}< t\leq t_{c3}, \\
0\,\text{mg/min} & \hbox{if}\,\, t_{c3}< t\leq t_{f}, \\
\end{cases}	
\end{equation} 
where $t_{c1}$, $t_{c2}$, and $t_{c3}$ are the three switching times. 
The trajectory $x(t)$ associated with Control 
\eqref{control:S7} is given by 
\begin{equation}
x(t)=\begin{cases}
\int\limits_0^{t}\exp(A(t-s))BU_{max}ds & \hbox{if}\,\, 0\leq t\leq t_{c1}, \\
\exp(A\,(t-t_{c1}))\,x(t_{c1})& \hbox{if}\,\, t_{c1}<t\leq t_{c2},	\\
\exp(A\,(t-t_{c2}))\,x(t_{c2})+\int\limits_{t_{c2}}^{t}\exp(A(t-s))BU_{max}ds 
& \hbox{if}\,\, t_{c2}<t\leq t_{c3},\\
\exp(A\,(t-t_{c3}))\,x(t_{c3})& \hbox{if}\,\, t_{c3}<t\leq t_{f}.	\\	
\end{cases}
\end{equation}
To compute the three switching times $t_{c1}$, $t_{c2}$, and $t_{c3}$ 
and the final time $t_f$, we have to solve the 
nonlinear system \eqref{non:linear:equation} with
\begin{equation}
\label{switching:time:tc1:tc2:S7}
0\leq t_{c1}\leq t_{c2}\leq t_{c3}\leq t_f.
\end{equation}
It turns out that System \eqref{non:linear:equation} subject to 
Condition \eqref{switching:time:tc1:tc2:S7} has no solution. 
Thus, Strategy~7 is also not feasible.

\begin{figure}[H]
\includegraphics[scale=0.38]{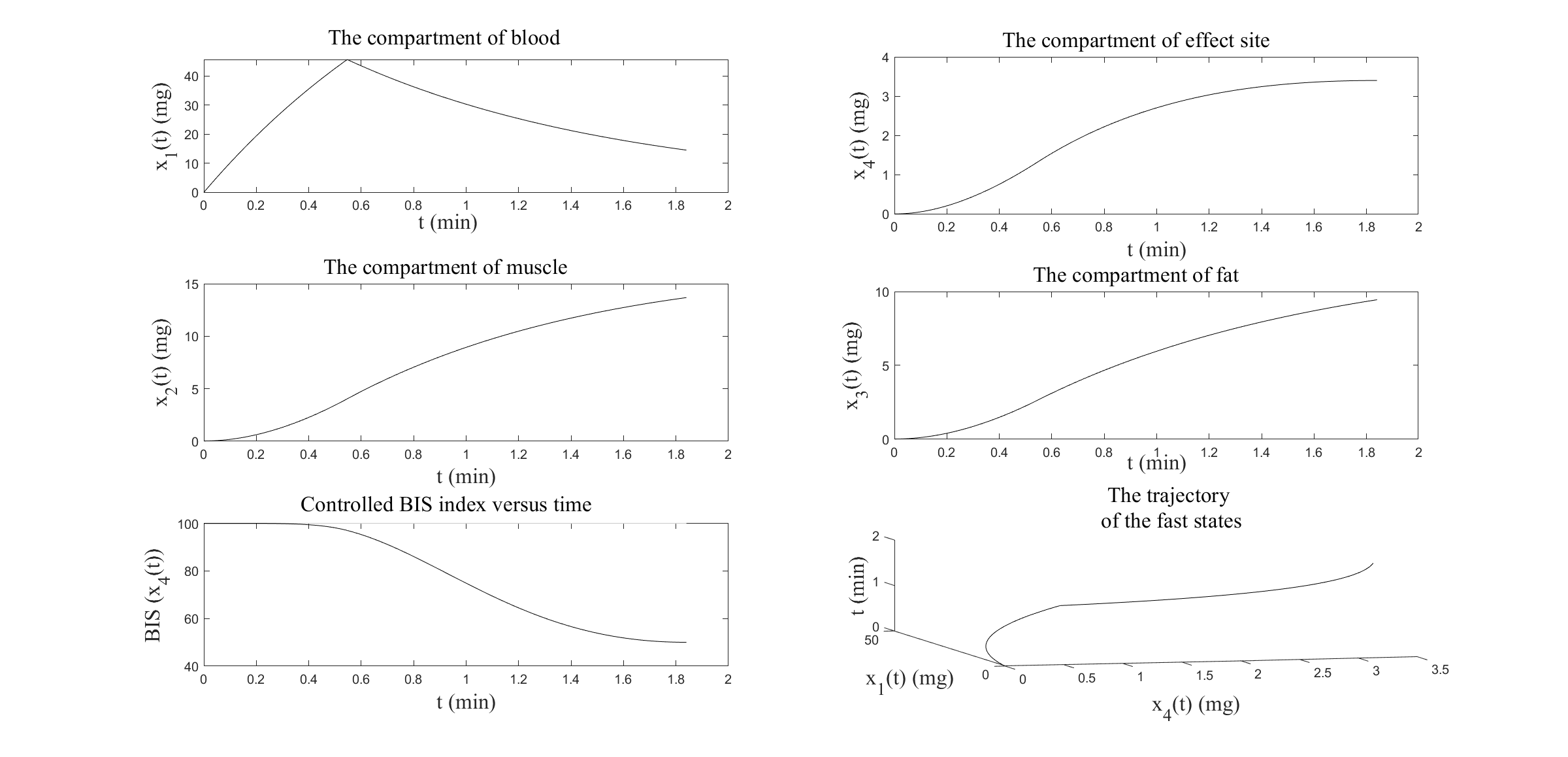}
\caption{The state trajectory, 
controlled BIS index, and trajectory of the fast 
states corresponding to the optimal control $u(t)$ 
of Figure~\ref{schema:optimal:control}, using the shooting method.}
\label{schema:Trajectory:shooting}
\end{figure}

\begin{figure}[H]
\includegraphics[scale=0.38]{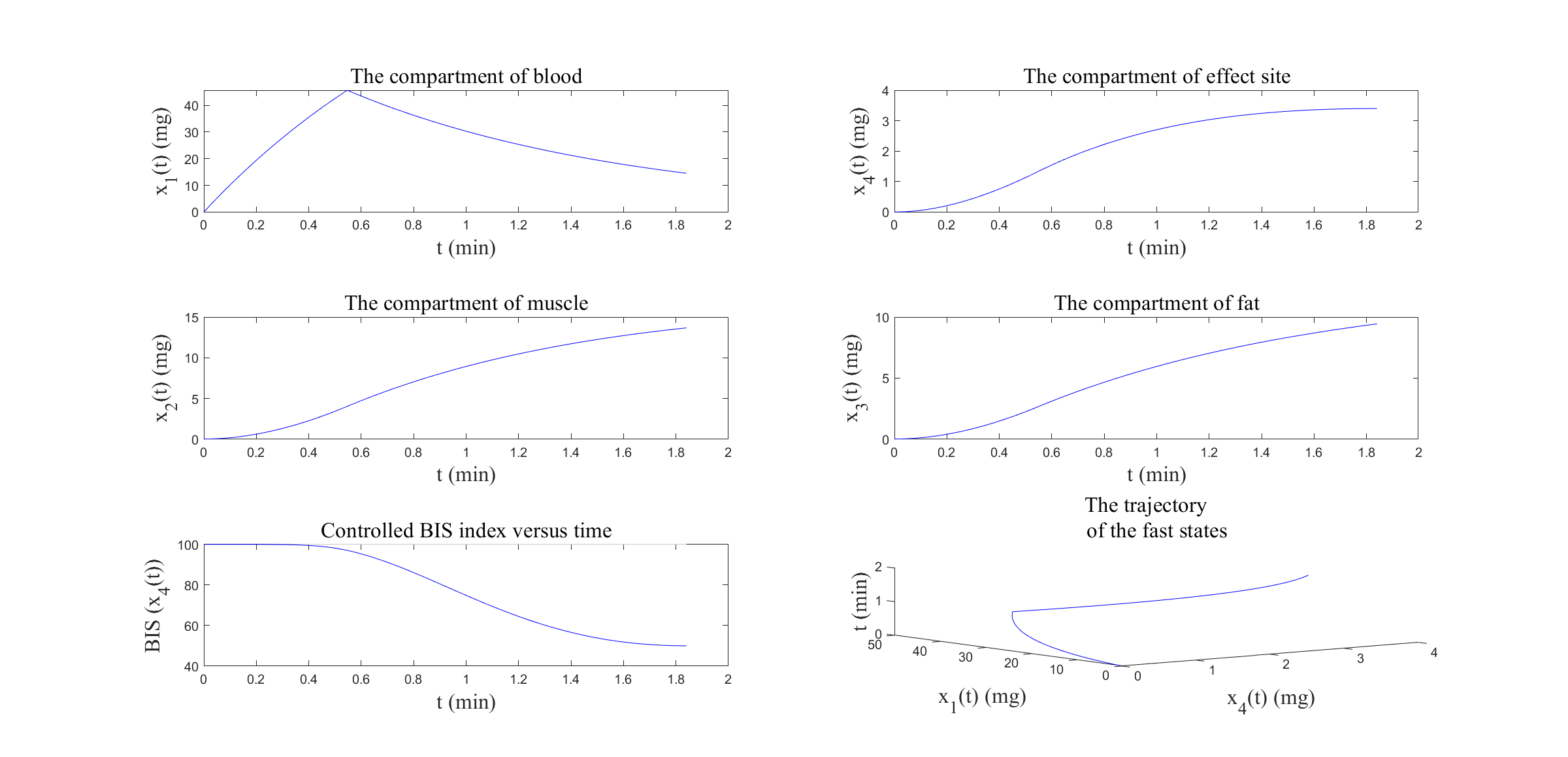}
\caption{The state trajectory, controlled BIS index, 
and trajectory of the fast states corresponding to the optimal control $u(t)$ 
of Figure~\ref{schema:optimal:control}, using the analytical method.}
\label{schema:Trajectory:analytical}
\end{figure}

\begin{figure}[H]
\includegraphics[scale=0.38]{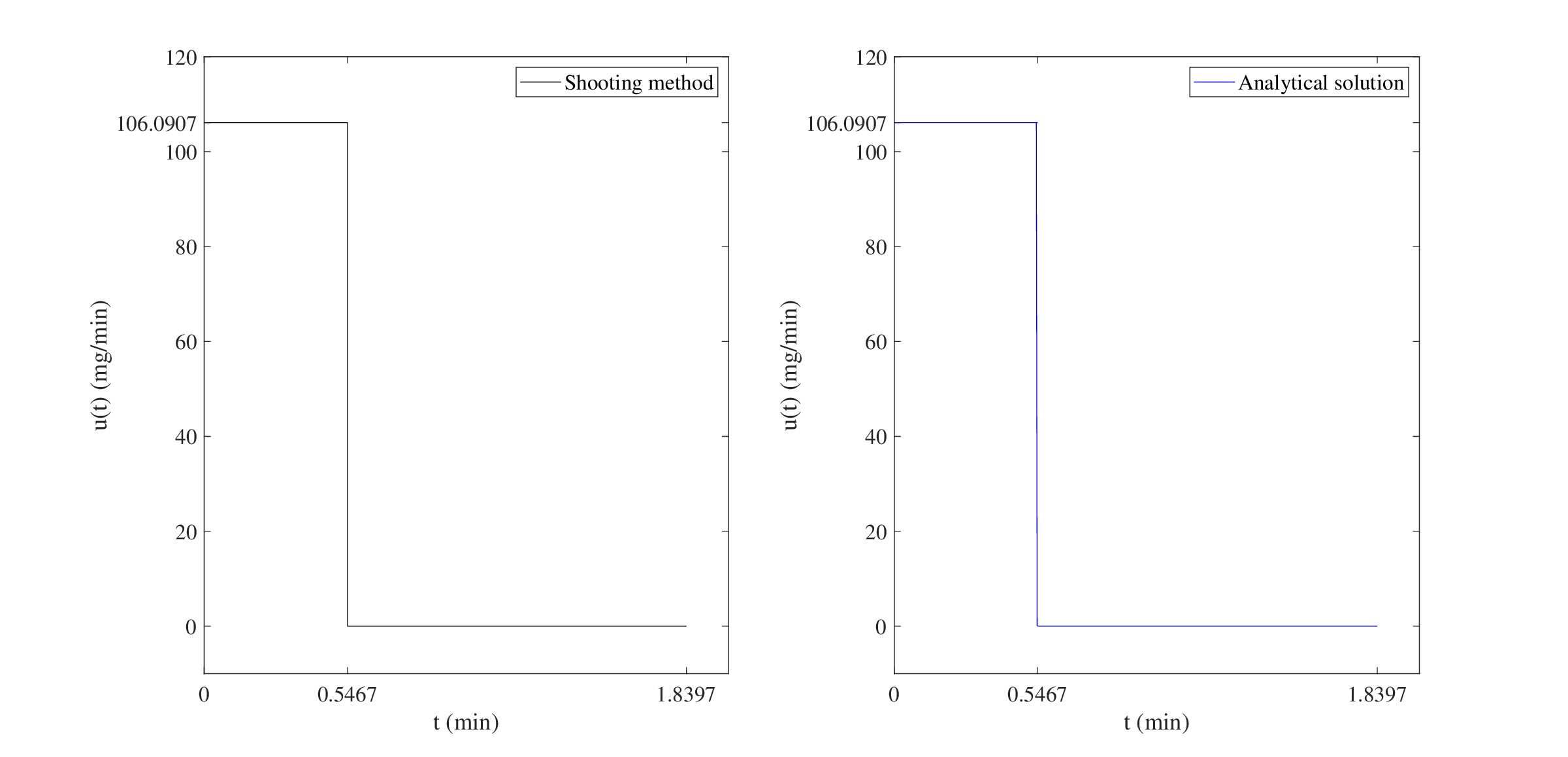}
\caption{The optimal continuous infusion rate $u(t)$ 
of the induction phase of anesthesia,
as obtained by the shooting and analytical methods.}
\label{schema:optimal:control}
\end{figure}

In Figures~\ref{schema:Trajectory:shooting} 
and \ref{schema:Trajectory:analytical}, we present the solutions 
of the linear system of differential equations \eqref{Example1:Time:Optimal} 
under the optimal control $u(t)$ illustrated in Figure~\ref{schema:optimal:control}, 
where the black curve corresponds to the one obtained by the shooting method, 
as explained in Section~\ref{Shooting:method}, while the blue curve corresponds 
to our analytical method, in the sense of Section~\ref{Analytical:solution}. In addition, 
for both figures, we show the controlled BIS Index, the trajectory of fast states 
corresponding to the optimal continuous infusion rate of the anesthetic $u(t)$, 
and the minimum time $t_f$ required to transition System \eqref{Example1:Time:Optimal} 
from the initial (wake) state
$$
x_{0}=(0,\,0,\,0,\,0)
$$ 
to the fast final (anesthetized) state
$$
x_{eF}=(14.518,\,3.4)
$$  
in the shortest possible time. The minimum time $t_f$ 
is equal to $t_f=1.8397\, \text{min} $ by the shooting method 
(black curve in Figure~\ref{schema:Trajectory:shooting}),
and it is equal to $t_f=1.8397\, \text{min}$ by the analytical 
method (blue curve in Figure~\ref{schema:Trajectory:analytical}).

By using the shooting method, the black curve in 
Figure~\ref{schema:optimal:control}
shows that the optimal continuous 
infusion rate of the induction phase of anesthesia $u(t)$ is equal to 
$106.0907$~mg/min until the switching time
$$
t_c=0.5467\, \text{min}.
$$
Then, it is equal to $0\,\text{mg/min}$ (stop-infusion) until the final time 
$$
t_f=1.8397\, \text{min},
$$

\vspace{-2pt}
\vspace{-8pt}


\vspace{-6pt}

By using the analytical method, 
the blue curve in Figure~\ref{schema:optimal:control}
shows that the optimal continuous infusion rate of 
the induction phase of anesthesia $u(t)$ 
is equal to $106.0907\,\text{mg/min}$ until the switching time
$$	
t_c=0.5467\, \text{min}.
$$
Then, it is equal to $0\,\text{mg/min}$ (stop-infusion) 
until the final time
$$
t_f=1.8397\, \text{min}.
$$
We conclude that both methods work well and give similar results.
However, in general, the shooting method does not always converge, 
depending on the initial conditions \eqref{eq:ic:psi}.
To obtain such initial values is not an easy task
since no theory is available to find them.
For this reason, the proposed analytical method
is logical, practical, and more suitable for real applications.


\section{Conclusions}
\label{sec:conc}

The approach proposed by the theory of optimal control is very effective. 
The shooting method was proposed by Zabi et al.~\cite{Said}, 
which is used to solve the time-optimal control problem and calculate the minimum time. 
However, this approach is based on Newton's method. The convergence of Newton's method 
depends on the initial conditions, being necessary to select an appropriate initial value 
so that the function is differentiable and the derivative does not vanish. This implies 
that the convergence of the shooting method is attached to the choice of the initial values.
Therefore, the difficulty of the shooting method is to find the initial conditions 
of the adjoint vectors. Here, the aim was to propose a different approach, which we call
``the analytical method'', that allows to solve the time-optimal control problem 
for the induction phase of anesthesia without such drawbacks. Our method is
guided by the selection of the optimal strategy,
without the need to choose initial values and study the convergence. 
We claim that our method can also be applied to other PK/PD models, 
in order to find the optimal time for the drug administration.

In the context of PK/PD modeling, the challenges associated with uncertainties 
in plant model parameters and controller gains for achieving robust stability 
and controller non-fragility are significant~\cite{Elsisi}. These challenges 
arise from factors like inter-individual variability, measurement errors, 
and the dynamic nature of patient characteristics and drug response. Further 
investigation is needed to understand and develop effective strategies to mitigate 
the impact of these uncertainties in anesthesia-related PK/PD models. This research 
can lead to the development of robust and non-fragile control techniques that 
enhance the stability and performance of anesthesia delivery systems. 
By addressing these challenges, we can improve the precision and safety of drug 
administration during anesthesia procedures, ultimately benefiting patient outcomes 
and healthcare practices. In this direction, the recent results of~\cite{Mohamed} 
may be useful. Moreover, we plan to investigate PK/PD fractional-order models, 
which is a subject under strong current research~\cite{MR4539083}. This is under 
investigation and will be addressed elsewhere.


\vspace{6pt} 

\authorcontributions{Conceptualization, M.A.Z., C.J.S., and D.F.M.T.; 
methodology, M.A.Z., C.J.S., and D.F.M.T.; 
software, M.A.Z.; 
validation, C.J.S. and D.F.M.T.; 
formal analysis, M.A.Z., C.J.S., and D.F.M.T.; 
investigation, M.A.Z., C.J.S., and D.F.M.T.; 
writing---original draft preparation, M.A.Z., C.J.S., and D.F.M.T.; 
writing---review and editing, M.A.Z., C.J.S. and D.F.M.T.; 
visualization, M.A.Z.; 
supervision, C.J.S. and D.F.M.T.; 
funding acquisition, M.A.Z., C.J.S., and D.F.M.T. 
All authors have read and agreed to the published version of the manuscript.}

\funding{This research was funded by the 
Portuguese Foundation for Science and Technology 
(FCT---Funda\c{c}\~{a}o para a Ci\^{e}ncia e a Tecnologia) 
through the R\&D Unit CIDMA, Grant Numbers UIDB/04106/2020 
and UIDP/04106/2020, and within the project ``Mathematical Modelling 
of Multi-scale Control Systems: Applications to Human Diseases'' (CoSysM3), 
Reference 2022.03091.PTDC, financially supported by national funds (OE) 
through FCT/MCTES.}

\institutionalreview{Not applicable.}

\informedconsent{Not applicable.}
	
\dataavailability{No new data were created or analyzed in this study. 
Data sharing is not applicable to this article. 
The numerical simulations of Section~\ref{Numerical:example} 
were implemented in \textsf{MATLAB} R2022a.	
The computer code is available from the authors upon request.}

\acknowledgments{The authors are grateful to four anonymous referees
for their constructive remarks and questions that helped to improve the paper.}

\conflictsofinterest{The authors declare no conflicts of interest. 
The funders had no role in the design of the study; in the collection, analyses, 
or interpretation of data; in the writing of the manuscript; 
or in the decision to publish the results.} 


\begin{adjustwidth}{-\extralength}{0cm}
	
\reftitle{References}


	
\PublishersNote{}

\end{adjustwidth}


\end{document}